\documentclass[11pt]{article}
\usepackage{amsmath}
\usepackage{amssymb,amsbsy,amsthm}
\usepackage[dvips]{graphicx}
\usepackage{caption}
\usepackage{subcaption}
\usepackage{pstricks}
\usepackage{pst-eps}
\usepackage{pst-node}
\usepackage{enumerate}
\usepackage{dsfont}
\usepackage[margin=1in]{geometry}
\usepackage{pdfsync}
\usepackage{hyperref}

\allowdisplaybreaks[1]

\newcommand{\p}{\mathbb{P}} % Probability P
\DeclareMathOperator{\poly}{poly} % Polynomial polys

\newcommand{\Beta}{\mathrm{Beta}} % Beta distribution
\newcommand{\Ga}{\mathrm{Gamma}} % Gamma distribution Ga
\newcommand{\GGa}{\mathrm{GGa}} % Generalized Gamma distribution GGa

\newtheorem{theorem}{Theorem}

\newtheorem{claim}[theorem]{Claim}

\begin{document}

\title{Beta-gamma tail asymptotics}
\author{
	Jim Pitman
	\thanks{University of California, Berkeley; \texttt{pitman@stat.berkeley.edu}.}
	\and 
	Mikl\'os Z.\ R\'acz
	\thanks{Microsoft Research; \texttt{miracz@microsoft.com}. Most of the work done while at University of California, Berkeley.}
}
\date{\today}

\maketitle

%%%%%%%%%%%%%%%%
%%% Abstract %%%
%%%%%%%%%%%%%%%%

\begin{abstract}
We compute the tail asymptotics of the product of a beta random variable and a generalized gamma random variable which are independent and have general parameters. 
A special case of these asymptotics were proved and used in a recent work of Bubeck, Mossel, and R\'acz 
in order to determine the tail asymptotics of the maximum degree of the preferential attachment tree. % started from an arbitrary initial tree. 
The proof presented here is simpler and highlights why these asymptotics hold. 
\end{abstract}

%%%%%%%%%%%%%%%%
%%% Document %%%
%%%%%%%%%%%%%%%%

%%%%%%%%%%%%%%%%%%%%%%%%%%%%%%%%%%%%%%%%%%%%
\section{Introduction} \label{sec:intro} %%%
%%%%%%%%%%%%%%%%%%%%%%%%%%%%%%%%%%%%%%%%%%%%

There has been a lot of recent interest in various urn schemes due to their appearance in many graph growth models 
(see, e.g.,~\cite{Mori:05,Ad-BeBrGo:10,Backhausz:11,PeRoRo:13,BeBoChSa:14,PeRoRo:14,BuMoRa:15,GoldschmidtHaas:15,PeRoRo:15+}). 
The limiting distributions arising in these urn schemes are often related to the beta and gamma distributions. 
Consequently, the computation of various statistics in random graph models often boils down to using algebraic properties of these distributions, 
commonly referred to as the beta-gamma algebra~\cite{Dufresne:98}.

The purpose of this note is to simplify and demystify a recent computation done in~\cite{BuMoRa:15} involving beta and generalized gamma random variables. 
Bubeck, Mossel, and R\'acz~\cite{BuMoRa:15} were interested in the influence of the seed graph in the preferential attachment model, 
which led them to study the tail asymptotics of the maximum degree of the preferential attachment tree. 
This, in turn, essentially reduces to computing the asymptotics of 
$\p \left( BZ > t \right)$ as $t \to \infty$, 
where $B \sim \Beta \left( a, b \right)$ and $Z \sim \GGa \left( a + b + 1, 2 \right)$ are independent random variables 
and $a$ and $b$ are positive integers; 
here $\Beta \left( a, b \right)$ denotes the beta distribution with positive parameters $a$ and $b$ 
(with density $\tfrac{1}{B \left( a, b \right)} x^{a-1} (1-x)^{b-1} \mathbf{1}_{\left\{ x \in [0,1]\right\}}$, 
where $B \left( a, b \right) = \tfrac{\Gamma \left( a \right) \Gamma \left( b \right)}{\Gamma \left( a+ b \right)}$ is the beta function), 
and $\GGa \left( c, p \right)$ denotes the generalized gamma distribution with density 
$\tfrac{p}{\Gamma \left( c/p \right)} x^{c-1} e^{-x^p} \mathbf{1}_{\left\{x > 0 \right\}}$ for $c,p > 0$. 
We refer to~\cite{BuMoRa:15} for details on these connections; see also~\cite{Jan06,PeRoRo:14}.

The computation in~\cite{BuMoRa:15} involves a few pages of alternating sums cancelling each other out in just the right way.   
Here, in contrast, we provide a short and simple proof of these asymptotics. 
The core calculation is only a few lines long, 
involving approximations at three points which are natural 
and which can be justified in a relatively straightforward manner. 
Moreover, the argument works for all positive values of the parameters $a$, $b$, $c$, and $p$. 
Throughout the paper we use standard asymptotic notation; for instance, 
$f\left( t \right) \sim g\left( t \right)$ as $t\to\infty$ 
if $\lim_{t\to\infty} f\left( t \right) / g \left( t \right) = 1$.

\begin{claim}\label{cl:tail}
Let $a$, $b$, $c$, and $p$ be positive, 
let $B \sim \Beta \left( a, b \right)$, 
and let $Z \sim \GGa \left( c, p \right)$, 
with $B$ and $Z$ independent. 
Then we have
\[
  \p \left( B Z > t \right) 
  \sim 
  \frac{\Gamma \left( a + b \right)}{\Gamma \left( c / p \right) \Gamma \left( a \right)} 
  p^{-b} t^{c - \left( b + 1 \right) p} 
  e^{-t^p}
  \tag*{as $t \to \infty$.}
\]
\end{claim}

There has been lots of work on understanding the distribution and tail asymptotics of products of random variables; 
see, e.g.,~\cite{SpringerThompson:70} for a paper from nearly half a century ago, 
and~\cite{HashorvaPakes:10} and references therein for recent developments. 
In particular, Claim~\ref{cl:tail} is a special case of~\cite[Theorem~4.1]{HashorvaPakes:10}, 
where the authors prove a general result for any product $BZ$ where $B \sim \Beta \left( a, b \right)$ and $Z$ has a law which is in the maximum domain of attraction of the Gumbel distribution. 
Due to the generality of their result their proof is fairly involved. 
We thus believe that the simple proof we present here is useful in highlighting why these asymptotics hold.

The product $BZ$ studied in Claim~\ref{cl:tail} has many nice properties, for instance it has moments of Gamma type~\cite{Janson:10}. 
When $p = 1$, $Z$ is a gamma random variable and $BZ$ has a so-called $G$ distribution, with its moments described by Meijer's $G$-function~\cite{Dufresne:10}. 
The case of $p = 2$ appears in many settings, including the preferential attachment model as mentioned above, and see also~\cite{Ad-BeBrGo:10,PeRoRo:15+} for connections to critical random graphs, random walks, and various random trees, including Aldous's Brownian continuum random tree (CRT). 
In a very interesting recent work, Pek\"oz, R\"ollin, and Ross~\cite{PeRoRo:15+} showed that generalized gamma random variables with $p$ being an integer greater than $2$ arise as limits in time inhomogeneous P\'olya-type urn schemes. 
It would be interesting to find connections to urn schemes and random graph models for general values of $p$.

%%%%%%%%%%%%%%%%%%%%%%%%%%%%%%%%%%%%%%%%%%%%%%%%%%%
\section{The core calculation} \label{sec:calc} %%%
%%%%%%%%%%%%%%%%%%%%%%%%%%%%%%%%%%%%%%%%%%%%%%%%%%%

In this section we prove Claim~\ref{cl:tail} modulo some approximations whose validity is justified later in Section~\ref{sec:approx}. 
Put $W = Z^p$, $z = w^{1/p}$, and $dz = \tfrac{1}{p} w^{1/p-1} dw$. 
The density of $W$ is thus 
\[
  f_W \left( w \right) = \frac{1}{\Gamma \left( c / p \right)} w^{c/p - 1} e^{-w} \mathbf{1}_{\left\{ w > 0 \right\}}.
\]
In other words, $W$ has a gamma distribution: $W \sim \Ga \left( c/p, 1 \right)$. 
We have 
\[
 \p \left( B Z > t \right) = \p \left( B^p Z^p > t^p \right) = \p \left( B^p W > t^p \right), 
\]
and so the claim is equivalent to 
\[
  \p \left( B^p W > w \right) 
  \sim 
  \frac{\Gamma \left( a + b \right)}{\Gamma \left( c/p \right) \Gamma \left( a \right)} 
  p^{-b} w^{c/p - b - 1} e^{-w}
  \tag*{as $w \to \infty$.}
\]
We can write 
\[
 \p \left( B^p W > w \right) = \p \left( W > w \right) \p \left( B^p W > w \, \middle| \, W > w \right). 
\]
The first factor has well-known asymptotics: 
$\p \left( W > w \right) \sim f_W \left( w \right)$ as $w \to \infty$~\cite[formula~6.5.32]{abramowitz1964handbook}. 
Also, it is well known that for any random variable $W$ which has a $\Ga \left( r, 1 \right)$ distribution, we have 
\[
 \left( W - w \, \middle| \, W > w \right) \stackrel{d}{\to} \mathcal{E} 
 \tag*{as $w \to \infty$,}
\]
where $\mathcal{E}$ is a standard exponential random variable. 
This convergence is rather strong, e.g., convergence of densities. 
So 
\begin{align}
\p \left( B^p W > w \, \middle| \, W > w \right) 
&= \p \left( B^p \left( w + W - w \right) > w \, \middle| \, W > w \right) \label{eq:approx1} \\ 
&\sim \p \left( B^p \left( w + \mathcal{E} \right) > w \right) 
= \p \left( B > \left( \frac{1}{1 + \mathcal{E} / w} \right)^{1/p} \right) \label{eq:approx2} \\
&\sim \p \left( B > 1 - \frac{\mathcal{E}}{pw} \right) 
= \p \left( 1 - B < \frac{\mathcal{E}}{pw} \right) \label{eq:approx3} \\
&\sim \int_0^\infty \frac{1}{B\left( a, b \right)} u^{b-1} e^{-pwu} du 
= \frac{1}{B \left( a, b \right)} \left( \frac{1}{pw} \right)^b \Gamma \left( b \right) \label{eq:approx4}
\end{align}
and the conclusion follows. 
The only thing that remains is to rigorously justify the three points where asymptotic equivalence was used in the line of reasoning above; 
see Section~\ref{sec:approx} for details.

%%%%%%%%%%%%%%%%%%%%%%%%%%%%%%%%%%%%%%%%%%%%%%%%%%%%%%%%%%%%%%
\section{Justifying the approximations} \label{sec:approx} %%%
%%%%%%%%%%%%%%%%%%%%%%%%%%%%%%%%%%%%%%%%%%%%%%%%%%%%%%%%%%%%%%

Here we justify why the expressions in~\eqref{eq:approx1},~\eqref{eq:approx2},~\eqref{eq:approx3}, and~\eqref{eq:approx4} are all asymptotically equivalent as $w \to \infty$. 

\textbf{Asymptotic equivalence of~\eqref{eq:approx3} and~\eqref{eq:approx4}.} 
By definition we have 
\[
  \p \left( 1 - B < \frac{\mathcal{E}}{pw} \right) 
  = \int_0^{\infty} \p \left( 1 - B < \frac{z}{pw} \right) e^{-z} dz 
  = \int_0^{pw} \p \left( 1 - B < \frac{z}{pw} \right) e^{-z} dz + e^{-pw}.
\]
Since $\int_0^\infty u^{b - 1} e^{-pwu} du = \Theta \left( w^{-b} \right)$ as $w \to \infty$,
we can neglect all terms that are $o \left( w^{-b} \right)$ as $w \to \infty$, 
and so we only have to deal with the integral term on the right hand side of the display above. 
Using that $1 - B \sim \Beta \left( b, a \right)$ and also the change of variables $z =  p w u$, we have
\begin{align*}
  \int_0^{pw} \p \left( 1 - B < \frac{z}{pw} \right) e^{-z} dz 
  &= \int_0^{pw} \int_0^{\frac{z}{pw}} \frac{1}{B \left( a, b \right)} x^{b-1} \left( 1 - x \right)^{a-1} dx e^{-z} dz \\ 
  &= \frac{pw}{B\left( a, b \right)} \int_0^1 \int_0^u x^{b-1} \left( 1 - x \right)^{a - 1} dx e^{-pwu} du. 
\end{align*}
The integral from $1/2$ to $1$ is negligible, since
\[
  pw \int_{1/2}^1 \frac{1}{B \left( a, b \right)} \int_0^u x^{b-1} \left( 1 - x \right)^{a-1} dx e^{-pwu} du 
  \leq 
  \int_{1/2}^1 \left( pw \right) e^{-pwu} du 
  \leq 
  e^{-pw/2}.
\]
For the integral from $0$ to $1/2$ we first drop the factor $\left( 1 - x \right)^{a-1}$ and we justify the validity of this later. 
Integration by parts then tells us that 
\begin{align*}
  \int_{0}^{1/2} \frac{u^b}{b} \left( pw \right) e^{-pwu} du 
  &= - \frac{1}{2^b b} e^{-pw/2} + \int_0^{1/2} u^{b-1} e^{-pwu} du \\
  &= \int_0^{\infty} u^{b-1} e^{-pwu} du - \frac{1}{2^b b} e^{-pw/2} - \int_{1/2}^{\infty} u^{b-1} e^{-pwu} du.
\end{align*}
Using that 
$\left( 1/2 + v \right)^{b-1} \leq 2 + \left( 2v \right)^{b-1}$ 
for all $b > 0$ and $v > 0$, we have that the integral from $1/2$ to $\infty$ in the display above is negligible: 
\begin{align*}
  \int_{1/2}^{\infty} u^{b-1} e^{-pwu} du 
  &= e^{-pw/2} \int_0^{\infty} \left( 1/2 + v \right)^{b-1} e^{-pwv} dv \\
  &\leq e^{-pw/2} \int_0^{\infty} \left( 2 + \left( 2 v \right)^{b-1} \right) e^{-pwv} dv 
  \leq \poly \left( w \right) e^{-pw/2},
\end{align*}
where $\poly \left( w \right)$ is some polynomial in $w$. 
We have thus shown that 
\[
  \frac{pw}{B\left( a, b \right)} \int_0^{1/2} \int_0^u x^{b-1}  dx e^{-pwu} du 
  \sim 
  \int_0^\infty \frac{1}{B\left( a, b \right)} u^{b-1} e^{-pwu} du
  \tag*{as $w \to \infty$.} 
\] 
Finally, to justify dropping the $\left( 1 - x \right)^{a-1}$ factor, 
note that 
$\left| 1 - \left( 1 - x \right)^{a-1} \right| \leq \max \left\{ a, 2 \right\} x$ 
for all $a > 0$ and $x \in \left( 0, 1/2 \right)$, 
and so 
\begin{multline*}
  \left| 
  \frac{pw}{B\left( a, b \right)} \int_0^{1/2} \int_0^u x^{b-1} \left( 1 - x \right)^{a - 1} dx e^{-pwu} du 
  - 
  \frac{pw}{B\left( a, b \right)} \int_0^{1/2} \int_0^u x^{b-1} dx e^{-pwu} du
  \right| \\
  \leq 
  \max \left\{ a, 2 \right\} 
  \frac{pw}{B\left( a, b \right)} \int_0^{1/2} \int_0^u x^{b} dx e^{-pwu} du
  = O \left( w^{- \left( b + 1 \right)} \right).
\end{multline*}

\textbf{Asymptotic equivalence of~\eqref{eq:approx2} and~\eqref{eq:approx3}.} 
We need to show that, as $w \to \infty$, 
\begin{equation}\label{eq:second}
  \left| \p \left( B > 1 - \frac{\mathcal{E}}{pw} \right) -  \p \left( B > \left( \frac{1}{1 + \mathcal{E} / w} \right)^{1/p} \right) \right| 
  = o \left( w^{-b} \right). 
\end{equation} 
Using that $1 - x/p \leq \left( 1 + x \right)^{-1/p} < 1 - x/p + \tfrac{p+1}{2p^2} x^2$ for all $x > 0$, 
we have that the left hand side of~\eqref{eq:second} is at most 
\[
  \p \left( B \in \left[ 1 - \frac{\mathcal{E}}{pw}, 1 - \frac{\mathcal{E}}{pw} + \frac{p+1}{2p^2} \frac{\mathcal{E}^2}{w^2} \right], \mathcal{E} \leq \frac{p}{p+1} w \right) 
  + 2 \p \left( \mathcal{E} > \frac{p}{p+1} w  \right). 
\]
We know that $\p \left( \mathcal{E} > \frac{p}{p+1} w \right) = e^{-\frac{p}{p+1}w}$, so we only have to deal with the first term in the sum above. 
Using the change of variables $z = pwu$ we have 
\begin{multline*}
  \p \left( B \in \left[ 1 - \frac{\mathcal{E}}{pw}, 1 - \frac{\mathcal{E}}{pw} + \frac{p+1}{2p^2} \frac{\mathcal{E}^2}{w^2} \right], \mathcal{E} \leq \frac{p}{p+1} w \right) \\
  \begin{aligned}
  &= \int_0^{\frac{p}{p+1} w} \int_{\frac{z}{pw}-\frac{p+1}{2p^2}\frac{z^2}{w^2}}^{\frac{z}{pw}} \frac{1}{B \left( a, b \right)} x^{b-1} \left( 1 - x \right)^{a-1} dx e^{-z} dz \\
  &= \frac{pw}{B \left( a, b \right)} \int_0^{\frac{1}{p+1}} \int_{u - \frac{p+1}{2} u^2}^{u} x^{b-1} \left( 1 - x \right)^{a-1} dx e^{-pwu} du \\
  &\leq \frac{\max \left\{ 1, \left( p / p + 1 \right)^{a-1} \right\}}{B \left( a, b \right)} \left( pw \right) 
  \int_0^{\frac{1}{p+1}} \int_{u - \frac{p+1}{2} u^2}^{u} x^{b-1} dx e^{-pwu} du.
  \end{aligned}
\end{multline*}
Using simple estimates we have 
$\int_{u - \frac{p+1}{2} u^2}^{u} x^{b-1} dx 
= \tfrac{1}{b} \left( u^b - \left( u - \tfrac{p+1}{2} u^2 \right)^b \right) 
\leq \frac{\left( p+1 \right) \left( b + 1 \right)}{2b} u^{b+1}$ 
for all $u \in \left( 0 , \tfrac{1}{p+1} \right)$ 
and so 
\begin{multline*}
  \left( pw \right) \int_0^{\frac{1}{p+1}} \int_{u - \frac{p+1}{2} u^2}^{u} x^{b-1} dx e^{-pwu} du 
  \leq \frac{\left( p+1 \right) \left( b + 1 \right)}{2b} \left( p w \right) \int_0^{\frac{1}{p+1}} u^{b+1} e^{-pwu} du \\
  \leq \frac{\left( p+1 \right) \left( b + 1 \right)}{2b} \left( p w \right) \int_0^{\infty} u^{b+1} e^{-pwu} du 
  =\frac{\left( p+1 \right) \left( b + 1 \right) \Gamma \left( b + 2 \right)}{2b p^{b+1}} w^{-\left( b + 1 \right)},
\end{multline*}
which concludes the proof of~\eqref{eq:second}. 

\textbf{Asymptotic equivalence of~\eqref{eq:approx1} and~\eqref{eq:approx2}.} 
We justify this via a direct calculation, though there might be a more elegant way to do this. 
To abbreviate notation, let $r := c/p$. By definition we have 
\begin{align*}
  \p \left( B^p \left( w + W - w \right) > w \, \middle| \, W > w \right) 
  &= \frac{\int_{w}^{\infty} \p \left( B^p \left( w + y - w \right) > w \right) \frac{1}{\Gamma \left( r \right)} y^{r-1} e^{-y} dy}{\frac{1}{\Gamma \left( r \right)} \int_w^{\infty} x^{r-1} e^{-x} dx} \\
  &= \frac{w^{r-1} e^{-w}}{\int_w^{\infty} x^{r-1} e^{-x} dx} 
  \int_0^{\infty} \p \left( B^p \left( w + z \right) > w \right) \left( 1 + \frac{z}{w} \right)^{r-1} e^{-z} dz. 
\end{align*}
For $r = 1$ this is exactly equal to 
$\p \left( B^p \left( w + \mathcal{E} \right) > w \right)$. 
Since
$\p \left( W > w \right) \sim f_W \left( w \right)$ as $w \to \infty$, 
the fraction in front of the integral in the display above goes to $1$ as $w \to \infty$, 
so what remains to show is that 
\[
  \int_{0}^{\infty} \p \left( B^p \left( w + z \right) > w \right) \left| \left( 1 + \frac{z}{w} \right)^{r-1} - 1 \right| e^{-z} dz 
  = o \left( \p \left( B^p \left( w + \mathcal{E} \right) > w \right) \right) 
  \tag*{as $w \to \infty$.} 
\] 
We partition the integral into two parts: 
from $0$ to $\sqrt{w}$, and from $\sqrt{w}$ to $\infty$. 
For the first term note that for $r \geq 1$ and $z \in \left( 0 , \sqrt{w} \right)$ 
we have 
$\left( 1 + z/w \right)^{r-1} - 1 \leq e^{r/\sqrt{w}} - 1$,
while for $r \in \left( 0 , 1 \right)$ and $z \in \left( 0, \sqrt{w} \right)$ 
we have 
$1 - \left( 1 + z/w \right)^{r-1} \leq 1 - \left( 1 + 1 / \sqrt{w} \right)^{-1} \leq 1 / \sqrt{w}$, 
and so 
\begin{multline*}
  \int_{0}^{\sqrt{w}} \p \left( B^p \left( w + z \right) > w \right) \left| \left( 1 + \frac{z}{w} \right)^{r-1} - 1 \right| e^{-z} dz \\
  \begin{aligned}
  &\leq 
  \max \left\{ \frac{1}{\sqrt{w}}, e^{\frac{r}{\sqrt{w}}} - 1 \right\} 
  \int_{0}^{\sqrt{w}} \p \left( B^p \left( w + z \right) > w \right) e^{-z} dz \\
  &\leq 
  \max \left\{ \frac{1}{\sqrt{w}}, e^{\frac{r}{\sqrt{w}}} - 1 \right\} 
  \p \left( B^p \left( w + \mathcal{E} \right) > w \right)
  = o \left( \p \left( B^p \left( w + \mathcal{E} \right) > w \right) \right). 
  \end{aligned}
\end{multline*} 
For the second term we can bound the factor 
$\p \left( B^p \left( w + z \right) > w \right)$ 
in the integral by $1$. 
For $r \in \left( 0, 1 \right)$ we have 
$1 - \left( 1 + z/w \right)^{r-1} \leq 1 - \left( 1 + z / w \right)^{-1} \leq z / w$
and so we get the upper bound of
\[
 \frac{1}{w} \int_{\sqrt{w}}^{\infty} z e^{-z} dz = \frac{1}{w} e^{-\sqrt{w}} \int_0^{\infty} \left( \sqrt{w} + u \right) e^{-u} du = \frac{\sqrt{w} + 1}{w} e^{-\sqrt{w}}. 
\]
For $r \geq 1$ we use the bound 
$\left( 1 + z/w \right)^{r-1} - 1 \leq e^{r z / w}$ 
to obtain the upper bound of 
\[
  \int_{\sqrt{w}}^{\infty} e^{\left( r / w - 1 \right) z} dz 
  = \frac{1}{1 - r/w} e^{\left( r/ w - 1 \right) \sqrt{w}} 
  \leq 2 e^{-\sqrt{w}/2},
\]
where we assumed that $w \geq 2r$.

%%%%%%%%%%%%%%%%%%%%%%%%
%%% Acknowledgements %%%
%%%%%%%%%%%%%%%%%%%%%%%%

\section*{Acknowledgements}

M.Z.R.\ gratefully acknowledges support from NSF grant DMS 1106999. 

%%%%%%%%%%%%%%%%%%
%%% References %%%
%%%%%%%%%%%%%%%%%%

\bibliographystyle{plain}
\bibliography{bib}

%%%%%%%%%%%%%%%%
%%% Appendix %%%
%%%%%%%%%%%%%%%%

% \newpage

% \appendix

\end{document}